\newcommand{\CLASSINPUTinnersidemargin}{17mm}
\documentclass[journal]{IEEEtran}
\ifCLASSINFOpdf \else \fi
\usepackage{amssymb,latexsym}
\usepackage{amsfonts}
\usepackage{amsthm,graphicx,epsfig}
\usepackage[cmex10]{amsmath}
\usepackage{cite}

\numberwithin{equation}{section}
\renewcommand{\theequation}{\thesection.\arabic{equation}}

\renewcommand{\baselinestretch}{1.1}

\begin{document}
\title{Effect of inter-sample spacing constraint on spectrum estimation with irregular sampling}
\author{Radhendushka Srivastava and~Debasis Sengupta~
\thanks{The authors are with the Applied Statistics Unit, Indian Statistical Institute, Kolkata,
700108, India (e-mail: radhe\_r@isical.ac.in;
sdebasis@isical.ac.in).}}

\maketitle

\IEEEpeerreviewmaketitle

\begin{abstract}
A practical constraint that comes in the way of spectrum estimation
of a continuous time stationary stochastic process is the minimum
separation between successively observed samples of the process.
When the underlying process is not band-limited, sampling at any
uniform rate leads to aliasing, while certain stochastic sampling
schemes, including Poisson process sampling, are rendered infeasible
by the constraint of minimum separation. It is shown in this paper
that, subject to this constraint, no point process sampling scheme
is alias-free for the class of all spectra. It turns out that point
process sampling under this constraint can be alias-free for
band-limited spectra. However, the usual construction of a
consistent spectrum estimator does not work in such a case.
Simulations indicate that a commonly used estimator, which is
consistent in the absence of this constraint, performs poorly when
the constraint is present. These results should help practitioners
in rationalizing their expectations from point process sampling as
far as spectrum estimation is concerned, and motivate researchers to
look for appropriate estimators of bandlimited spectra.
\end{abstract}


\section{Introduction}\label{S1}
Estimation of the power spectral density of a continuous time mean
square continuous stationary stochastic process is a classical
problem \cite{Kay}. Estimates are usually based on finitely many
observed points. If the spectral density of the process is
bandlimited (i.e., compactly supported), then one can estimate it
consistently by using uniformly spaced samples, provided the
sampling rate is faster than the Nyquist rate \cite{Shannon}. On the
other hand, if the spectral density is not bandlimited, then uniform
sampling at any sampling rate leads to aliasing, i.e., there is a
class of continuous processes whose spectral densities are
indistinguishable from the sampled process. Thus, one can not
estimate the spectral density of the original process consistently
on the basis of uniformly spaced samples \cite{Silverman}.

In such a situation, non-uniform or irregular sampling schemes have
been explored. Silverman and Shapiro \cite{Silverman} introduced a
notion of alias-free sampling.
Beutler \cite{Beutler} further formalized this definition of
alias-free sampling for different classes of power spectra. Masry
\cite{Masrynew} gave another definition of alias-free sampling. For
each sampling scheme that is alias-free according to this
definition, he provided a corresponding estimator of the spectral
density, which would be consistent under certain conditions. Poisson
sampling (i.e., sampling at the arrival times of a homogeneous
Poisson process) turns out to be alias-free for the class of all
spectra, according to both the definitions.

The existence of consistent estimators in the non-bandlimited case
makes non-uniform sampling schemes, in particular Poisson process
sampling, very attractive. However, consistency is only a large
sample property of an estimator. Srivastava and Sengupta \cite{our}
showed that, if one has the ability to sample the process
arbitrarily fast, then one can consistently estimate a
non-bandlimited spectral density through uniformly spaced samples
also, provided the sampling rate goes to infinity at a suitable rate
as the sample size goes to infinity. By comparing the smoothed
periodogram estimator with the corresponding estimator based on the
Poisson process sampling, they found that, under certain regularity
conditions, the rates of convergence for the two estimators are
comparable and the constants associated with the rates of
convergence have a trade-off in terms of bias and variance. Thus,
the existence of a consistent estimator is not an exclusive
advantage of point process sampling.

An attractive property of the spectrum estimator based on Poisson
sampling \cite{Masrypoisson} is that the estimator is consistent for
any average sampling rate. This property implies that one can
estimate bandlimited spectra consistently, even if the sampling is
done at a sub-Nyquist average rate. Such estimators show that
deficiencies in sampling rate can be made up by sample size,
provided one is prepared to sample at irregular intervals. This fact
gives rise to the hope that even when there is a constraint on the
sampling rate, one can judiciously use non-uniform sampling to
consistently estimate spectra with much larger bandwidth than what
can be achieved through uniform sampling.

It is important to note that a small average sampling rate does not
mean that {\it any} two successive samples are far apart. In the
case of Poisson process sampling with any average sampling rate, it
can be seen that as the sample size goes to infinity, there would be
a large number of pairs of consecutive samples which are nearer to
each other than any specified threshold. Thus, in order to use
Poisson process sampling with any average sampling rate, sometimes
one has to sample the process very fast. Many other non-uniform and
alias-free sampling schemes also have this requirement. All these
schemes become infeasible if there is a hard limit on the minimum
separation between successive samples. Such a constraint can arise
because of technological limits as well as economic considerations.

Books on sampling \cite{Higg,BenFer} give a clear picture of the
limitations of uniform sampling in respect of a constraint on the
minimum separation between successive samples. However, suitability
of non-uniform sampling schemes in the presence of this constraint
has not been studied so far \cite{Marvasti}.

In this paper, we consider the problem of consistent estimation of
the power spectral density of a stationary stochastic process
through non-uniform sampling, under a constraint on the minimum
separation between successive samples. In Section~\ref{S2}, we
describe the underlying set-up and discuss the notions of alias-free
sampling provided by Shapiro and Silverman \cite{Silverman} as well
as by Masry \cite{Masrynew}. In section~\ref{S3}, we consider the
class of all power spectra, and show that under the above
constraint, no stationary point process sampling scheme is
alias-free for this class. Subsequently, we study the possibility of
alias-free sampling for estimation of spectra that are known to be
confined to a certain bandwidth. In section~\ref{S4}, we discuss the
difficulties of obtaining a consistent estimator of the power
spectrum even when it is known to be bandlimited. In
Section~\ref{S5}, we report the results of a simulation study of the
performance of a commonly used estimator based on Poisson process
sampling, in the presence of the above constraint. We summarize the
findings and provide some concluding remarks in Section~\ref{S6}.

\section{Notions of Alias-free Sampling}\label{S2}

Let $X=\{X(t),~-\infty<t<\infty\}$ be a real, mean square continuous
and wide sense stationary stochastic process with mean zero,
covariance function $C(\cdot)$ and spectral distribution function
$\Phi(\cdot)$. If $\Phi(\cdot)$ has a density, we denote it by
$\phi(\cdot)$. Let $\tau=\{t_n,~n=\ldots,-2,-1,0,1,2,\ldots\}$ be a
sequence of real-valued, stochastic sampling times.

Shapiro and Silverman's \cite{Silverman} notion of alias free
sampling, which was further formalized by Beutler \cite{Beutler}, is
based on the following assumptions about the sampling process.

\bigskip
{Assumption \bf A1.} The process $\tau$ is independent of $X$.

{Assumption \bf A2.} The sequence of sampling times $\tau$
constitutes a stationary point process, such that the probability
distribution of $(t_{m+n}-t_m)$ does not depend on $m$.

\bigskip
Under a sampling scheme that satisfies properties {\bf A1} and {\bf
A2}, the sampled process $X_s=\{X(t),~t\in\tau\}$ is wide sense
stationary. Denote the covariance sequence of the sampled process
$X_s$ by $R=\{\ldots,r(-2),r(-1),r(0),r(1),r(2),\ldots\}$, where
$$r(n)=E[X(t_{m+n})X(t_m)],\ \mbox{for $m$, $n$ integers},$$
and the expectation is taken without conditioning on the sampling
times. Beutler's definition of alias-free sampling, based on Shapiro
and Silverman's earlier idea, is as follows.

\bigskip
\noindent {\bf Definition 1.} The sampling process $\tau$ satisfying
assumptions {\bf A1} and {\bf A2} is alias-free relative to the
class of spectra $\mathbb{S}$ if no two random processes with
different spectra belonging to $\mathbb{S}$ yield the same
covariance sequence ($R$) of the sampled process.

\bigskip
Shapiro and Silverman \cite{Silverman} had considered the special
case where the sampling times constitute a renewal process, and
$\mathbb{S}$ is the class of {\it all} spectra with integrable and
square integrable densities. They referred to this scheme as
additive random sampling, and showed that it is alias-free, provided
the characteristic function of the inter-arrival distribution takes
no value more than once on the real line. In particular, Poisson
process (a renewal process having exponentially distributed
inter-arrival times) sampling scheme is alias-free for the class of
spectra $\mathbb{S}$.

The above definition has the drawback that it does not make use of
the information contained in the sampling times. If one wishes to
reconstruct $\phi(\cdot)$ using a sampling scheme that is alias-free
according to the above definition, then that would be done on the
basis of the sequence $R$ only. Beutler \cite{Beutler} gave a
procedure for this reconstruction, and indicated that this procedure
may be used to estimate $\phi(\cdot)$ from estimates of $R$.
However, Masry \cite{Masrynew} pointed out that the above definition
does not lead to a spectrum estimator that is provably consistent.

From all these considerations, this approach appears to be rather
restrictive. In practice, one would expect to use the information
contained not only in the sampled values, but also in the sampling
times, in order to estimate the power spectral density. In order to
take into account the sampling times, Masry \cite{Masrynew} gave an
alternative definition of alias-free sampling, while making
Assumption {\bf A1} and the following additional assumption about
the sampling process.

\bigskip
{Assumption \bf A3.} The process $\tau$ constitutes a stationary
orderly second-order point process on the real line.

\bigskip
Let $\beta$ be the mean intensity and $\mu_{c}$ be the reduced
covariance measure of the process $\tau$, and $\mathbb{B}$ be the
Borel $\sigma$-field on the real line. Consider the compound process
$\{Z(B),~B\in\mathbb{B}\}$ defined by
$$Z(B)=\sum_{t_i\in B}X(t_i).$$
The process $Z=\{Z(B),~B\in\mathbb{B}\}$ is second order stationary
(i.e., the first and second moments of $Z(B+t)$, for any real number
$t$, does not depend of $t$). Let $\mu_z$ be the covariance measure
of the process $Z$. It can be shown that this measure is given by
\begin{equation}
\mu_z(B)=\int_{B}C(u)[\beta^2 du+\mu_c(du)].
\label{muz}
\end{equation}
Masry's notion of alias-free sampling is as follows.

\bigskip
\noindent {\bf Definition 2.} The sampling process $\tau$ satisfying
assumptions {\bf A1} and {\bf A3} is alias-free relative to the
class of spectra $\mathbb{S}$ if no two random processes with
different spectra belonging to $\mathbb{S}$ yield the same
covariance measure ($\mu_z$) of the compound process.

\bigskip
Note that this definition makes use of the information contained in
the sampling times, as the covariance measure $\mu_z$ involves the
mean intensity $\beta$ as well as the reduced covariance measure
$\mu_{c}$ of the sampling process. It has been shown that, according
to Definition~2, Poisson process sampling is alias-free for the
class of all spectra having integrable and square integrable
densities \cite{Masrynew}.

\section{Sampling under Constraint}\label{S3}

As mentioned in Section~\ref{S1}, the focus of the present work is
on a sampling process $\tau$ which satisfies the following
constraint.

\bigskip
{Assumption \bf A4.} The time separation between two successive
sample points is at least $d$ (i.e., $t_{n+1}-t_n\ge d$ for any
index $n$) for some fixed $d>0$.

\bigskip
In this section, we investigate whether a sampling scheme satisfying
this constraint can be alias-free.

\subsection{General spectra}
We present some negative results in the case when $\mathbb{S}$ is
the class of {\it all} spectra -- bandlimited or otherwise.

\bigskip\noindent
{\bf Theorem 1.} No sampling point process satisfying Assumptions
{\bf A1}, {\bf A2} and {\bf A4} is alias-free according to
Definition~1, for the class of all spectra.

\bigskip\noindent
{\bf Theorem 2.} No sampling point process satisfying Assumptions
{\bf A1}, {\bf A3} and {\bf A4} is alias-free according to
Definition~2, for the class of all spectra.

\bigskip
We prove these theorems in the appendix by constructing
counter-examples, based on the following class of power spectral
densities.
\begin{equation}
\mathbb{A}=\left\{\phi(\cdot):\int_{-\infty}^\infty\phi(\lambda)e^{it\lambda}d\lambda=0\mbox{
for }|t|>d.\right\} \label{classA}
\end{equation}
The members of this class correspond to covariance functions
supported over the interval $[-d,d]$. A member of this class is the
power spectral density defined by
$$\phi_a(\lambda)=\frac{1}{\pi a}\frac{1-\cos(a\lambda)}{\lambda^2},$$
for any arbitrary positive $a\in(0,d]$. This density corresponds to
the covariance function
$$C_a(t)=\begin{cases}
1-\frac{|t|}{a}&\mbox{for $|t|\le a$,}\\ 0&\mbox{for $|t|>a$.}\\
\end{cases}$$ %
Some other members of $\mathbb{A}$ can be constructed by convolving
$\phi_a(\cdot)$ with an arbitrary power spectral density. We show in
the appendix that if $X_1$, $X_2$ and $X_3$ are independent mean
square continuous stochastic processes such that $X_2$ and $X_3$
have spectra in $\mathbb{A}$ and have the same variance, then the
spectra of $X_1+X_2$ and $X_1+X_3$ cannot be distinguished from the
sequence $R$ or the measure $\mu_z$, leading to aliasing according
to Definitions~1 and~2.

One can easily construct two {\it integrable and square integrable}
power spectral densities that are indistinguishable from $R$ or
$\mu_z$. Therefore, the statements of Theorems~1 and~2 also hold in
respect of all spectra having integrable and square integrable
densities (rather than all spectra). Thus, the alias-free property
of Poisson process sampling mentioned in Section~\ref{S2} become
inapplicable, once the inter-sample spacings are adjusted in
accordance with the constraint~{\bf A4}.

These two theorems show that, under the constraint of a minimum
inter-sample spacing, any point process sampling scheme would be
inadequate for the identification of a completely unrestricted power
spectral density -- according to the existing notions of alias-free
sampling. If the power spectral density of the original continuous
time process is not identifiable from the sequence $R$ or the
covariance measure $\mu_z$, then one cannot expect to consistently
estimate the power spectral density on the basis of estimates of
either of these.

Note that the Assumption {\bf A4} comes from a practical
consideration, and it is difficult to think of an implementable
sampling scheme that would not require it (i.e., a scheme that can
have arbitrarily closely spaced samples).

It is well known that estimators based on uniformly spaced samples,
irrespective of the sampling rate, also suffer from the limitation
of non-identifiability. In fact, it is this limitation of uniform
sampling that has been historically used as one of the major
arguments in favor of non-uniform sampling schemes. The above
theorems show that the same difficulty applies to practical
non-uniform sampling schemes as well.

\subsection{Bandlimited Spectra}

In the case of uniform sampling, it is well known that a bandlimited
process would not lead to aliasing provided that the sampling is
done at the Nyquist rate or faster. On the other hand, uniform
sampling at any fixed rate would be free from the problem of
aliasing if and only if the spectrum of the continuous time process
is known to be confined to an appropriate band. This fact, together
with the limitation of point process sampling in the case of
non-bandlimited spectra, gives rise to the question: Can point
process sampling under Assumption {\bf A4} be alias-free for the
class of bandlimited spectra? If so, it would be interesting to
compare the maximum allowable spectral bandwidths for alias-free
sampling, arising from uniform and point process sampling schemes
under Assumption {\bf A4}.

It turns out that alias-free sampling under Assumption {\bf A4} is
possible for an important class of stochastic sampling schemes,
namely, renewal process sampling. This is a special case of point
process sampling, which has received much attention from researchers
\cite{Silverman,Beutler,Masrynew,Tar}. Poisson process sampling is a
further special case of renewal process sampling. However, it is an
ideal sampling scheme, in contrast with implementable schemes that
would require Assumption {\bf A4}.

\begin{figure}[!t]
\includegraphics[width=3.5in,height=2.5in]{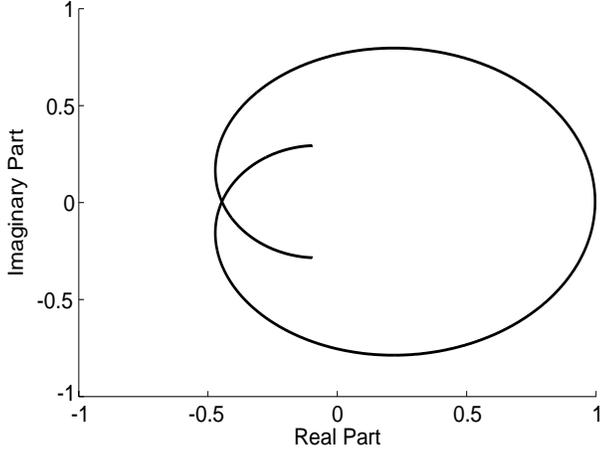}
\caption{Graph of $f'(\lambda)$ on the complex plane for the left
truncated exponential distribution, with mean $2d$ and truncation
point $d$, for $-\pi/d<\lambda\le\pi/d$.}
\end{figure}

The benchmark for the present study would be the fastest possible
rate of uniform sampling under Assumption {\bf A4}, which is $1/d$.
Note that uniform sampling at this rate is alias-free for the class
of spectra supported on $[-\pi/d,\pi/d]$.

First, we present a general result that would be useful in answering
the foregoing question, as far as Definition 1 of alias-free
sampling is concerned.

\bigskip\noindent
{\bf Theorem 3.} A renewal process sampling scheme satisfying
Assumptions {\bf A1} and {\bf A2}, and having characteristic
function of the inter-sample spacing denoted by $f'$, is alias-free
relative to a class of spectra supported on the closed and finite
interval $I$ according to Definition~1 if and only if the graph of
$f'(\lambda)$ on the complex plane, for $\lambda\in I$, does not
divide the complex plane.

\bigskip\noindent
Theorem~3 relates the alias-free property of a renewal process
sampling scheme to the geometry of the characteristic function of
the inter-sample spacing. It may be noted that the distribution of
$d+X$, where $d$ is fixed and $X$ has the gamma distribution with
any combination of parameters, does not satisfy the necessary and
sufficient condition given in Theorem 3 for $I=[-\pi/d,\pi/d]$. It
follows that the corresponding renewal process sampling schemes,
including the case of inter-sample spacing having a left-truncated
exponential distribution, are not alias-free according to Definition
1, relative to a class of spectra limited to the band
$[-\pi/d,\pi/d]$. The graph of $f'(\lambda)$ for the left-truncated
exponential distribution with mean $2d$ and truncation point $d$,
for $-\pi/d\le \lambda\le\pi/d$, is shown in Figure 1. For such
sampling schemes, aliasing can be avoided only if the continuous
time process is confined to a bandwidth that is even smaller than
the maximum allowable bandwidth in the case of uniform sampling.

However, there are some other renewal process sampling schemes that
satisfy Assumption {\bf A4} and are alias-free for the class of
spectra limited to a band larger than $[-\pi/d,\pi/d]$, as the next
theorem shows.

\begin{figure}[!t]
\includegraphics[width=3.5in,height=2.5in]{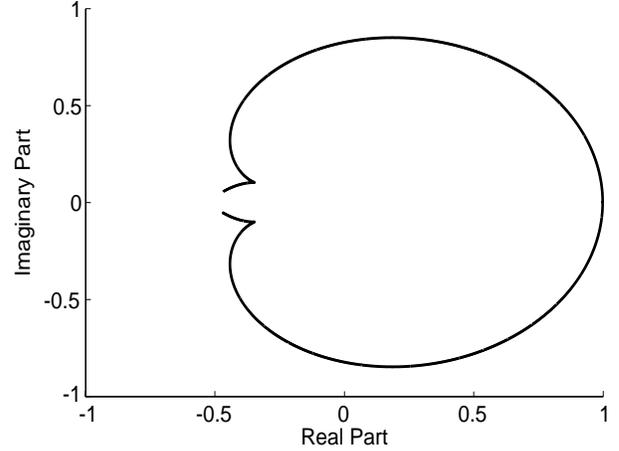}
\caption{Graph of $f'(\lambda)$ on the complex plane for the example
given in the proof of Theorem~4, for
$-1.1\pi/d<\lambda\le1.1\pi/d$.}
\end{figure}

\bigskip\noindent
{\bf Theorem 4.} There exists a closed and finite interval $I$,
which contains the interval $[-\pi/d,\pi/d]$, and a renewal process
sampling scheme satisfying Assumptions {\bf A1}, {\bf A3} and {\bf
A4} which is alias-free relative to the class of spectra supported
on $I$, according to Definition~1.

\bigskip
The proof of Theorem~4, given in the appendix, invokes an example,
for which $I$ is more than 10\% larger than the interval
$[-\pi/d,\pi/d]$, while the average inter-sample spacing is about
35\% more than the minimum allowable spacing ($d$). The graph of
$f'(\lambda)$ for this inter-sample spacing distribution, for
$-1.1\pi/d\le \lambda\le1.1\pi/d$, is shown in Figure 2.

\bigskip
We now turn to Definition 2. Since this notion of alias-free
sampling is weaker than that of Definition 1, one can expect a
stronger result.

\bigskip\noindent
{\bf Theorem 5.} Any renewal process sampling scheme, satisfying
Assumptions {\bf A1}, {\bf A3} and {\bf A4} and the further
assumption that the inter-sample spacing distribution has a density
that is positive over a semi-infinite interval, is alias-free
according to Definition~2, for the class of spectra limited to the
band $\left[-\lambda_0,\lambda_0\right]$ for every finite
$\lambda_0>0$.

\bigskip Theorem 5 shows that, under the constraint of a minimum allowable
separation between successive samples, renewal process sampling is
alias-free (according to Definition~2) for a {\it wider} range of
power spectra than uniform sampling. It is interesting to note that
sampling schemes following the assumptions of Theorem 5 are
alias-free according to Definition~2 when the spectral density of
the underlying continuous-time process is known to be confined to
any finite bandwidth (no matter how large), but according to Theorem
2, these are not alias-free when the process is non-bandlimited.

It transpires from the foregoing discussion that there are
contrasting scopes of alias-free renewal process sampling under the
constraint of a minimum allowable inter-sample spacing, according to
Definitions~1 and~2. The limited scope of alias-free sampling in the
case of Definition 1 stems from the fact that, under that notion,
one aims to identify spectra solely from the sequence $R$, which is
rather restrictive.

\section{Difficulties in estimation of bandlimited Spectrum}\label{S4}

Consider a class of spectra having density supported on the closed
and finite interval $I$. Given a sampling scheme that is alias-free
relative to this class according to Definition 1, one would look for
an estimate of the power spectral density based on estimated values
of the sequence $R$. Beutler \cite{Beutler} outlined a method of
estimation based on the representation
\begin{equation}\label{rn}
\Phi(\lambda_0)=\lim_{n\rightarrow\infty}\sum_{k=1}^{n}c_{kn}r(n),
\end{equation}
for a continuity point $\lambda_0$ of the spectral distribution
function, where $c_{kn}$, $k=1,\ldots,n$ are the coefficients such
that the uniform convergence of the sequence of partial sums
\begin{equation}\label{in}
\lim_{n\rightarrow\infty}\sum_{k=1}^{n}c_{kn}f'_{k}(\lambda)=1_{\left\{(-\infty,\lambda_0)\cap
I\right\}}(\lambda)
\end{equation}
happens everywhere except perhaps at $\lambda_0$. Here, for
$k=1,2,\ldots$, $f'_{k}(\lambda)$ is the characteristic function of
the distribution of the spacing between $k$ successive samples. Note
that the representation (\ref{in}) is possible whenever the sampling
scheme is alias-free according to Definition~1\cite{Beutler}.

One can estimate the spectral distribution $\Phi(\cdot)$ by plugging
in the estimate of $r(n)$ in (\ref{rn}), and can subsequently obtain
an estimator of the spectral density $\phi(\cdot)$. The coefficients
$c_{kn}$ defined by (\ref{in}), however, are attributes of the
sampling scheme, and these have to be computed theoretically. For
the example constructed in the proof of Theorem~4, the
characteristic function happens to be
$$f'_k(\lambda)=[0.68e^{i\lambda d}+0.32e^{i2.1\lambda d}]^{k}.$$
The coefficients $c_{1n},\ldots,c_{nn}$ for any fixed $n$ can be
obtained numerically, either directly from (\ref{in}) or by using a
Gram-Schmidt orthogonalization of the characteristic functions as
suggested in \cite{Masrynew}. There is no closed form solution. In
any case, the consistency of the plug-in estimator based on
(\ref{rn}) has not been proved.

Theorem~5 gives us a reason to look for estimators of $\phi(\cdot)$
based on constrained sampling schemes that are alias-free for
bandlimited processes according to Definition~2. Such an estimator
would be based on the measure $\mu_z$. It can be shown that the
power spectral density $\phi(\cdot)$ is related to the
characteristic functions $\phi_z$ and $\phi_c$ of the measures
$\mu_z$ and $\mu_c$ defined as
\begin{eqnarray*}
\phi_z(\lambda)&=&\frac{1}{2\pi}\int_{-\infty}^{\infty}e^{-iu\lambda}\mu_z(du),\\
\phi_c(\lambda)&=&\frac{1}{2\pi}\int_{-\infty}^{\infty}e^{-iu\lambda}\mu_c(du),
\end{eqnarray*}
respectively, through the integral equation
\begin{equation}
\phi_z(\lambda)=\beta^{2}\phi(\lambda)
+\int_{-\infty}^{\infty}\phi(\lambda-\omega)\phi_c(\omega)d\omega,\label{integral}
\end{equation}
provided that the measure $\mu_c$ is totally finite. Masry
\cite{Masrynew} and Brillinger \cite{Brillinger} gave a unique and
explicit solution to (\ref{integral}) under the following additional
assumption about the sampling process~$\tau$.

\bigskip\noindent
Assumption {\bf B1.} The reduced covariance measure $\mu_c$ has the
density $f_c(\cdot)$ satisfying the conditions\\
\begin{enumerate}
\item[(i)] $f_c(u)+\beta^{2}>0$ for all
$u\in(-\infty,\infty)$;\\[.5ex]
\item[(ii)] $\frac{f_c(u)}{f_c(u)+\beta^{2}}$  is integrable.\\
\end{enumerate}
Masry \cite{Masrynew} also showed that an empirical version of this solution is a
consistent estimator of $\phi(\cdot)$.

Note that, in the case of renewal processes,
$(f_c(\cdot)+\beta^2)/\beta$ is the renewal density. Under
Assumption {\bf A4}, this density would be necessarily zero over the
interval $[-d,d]$. Therefore, Assumption {\bf B1} is violated, and
the solution to (\ref{integral}), given in
\cite{Masrynew,Brillinger}, is not applicable to the present
situation. No explicit solution is available in general.
Consequently, there is no scope of using Masry's plug-in estimator.

It emerges from this discussion that all the estimators of a
bandlimited power spectral density, which have been proposed so far
on the basis of renewal process sampling, are either inapplicable in
the present context, or cannot be shown to be consistent.

In summary, even though renewal process sampling, subject to the
constraint of a minimum inter-sample spacing, is alias-free for
bandlimited power spectral densities, there is no estimator in the
literature that is known to be consistent. One has to look for new
estimators that may be appropriate in this situation.

\section{Simulation}\label{S5}
The foregoing discussion leads us to an interesting question: How
would the {\it available} estimators perform under the constraint of
a minimum inter-sample spacing? The performance of various
estimators based on uniform sampling have been studied both
theoretically and empirically, and their limitations arising from
aliasing have been exposed. In this section, we consider the
performance of a well-known estimator based on non-uniform sampling,
under the constraint of a minimum inter-sample spacing.

\begin{figure*}[!t]
\noindent
\includegraphics[width=7in,height=5.5in]{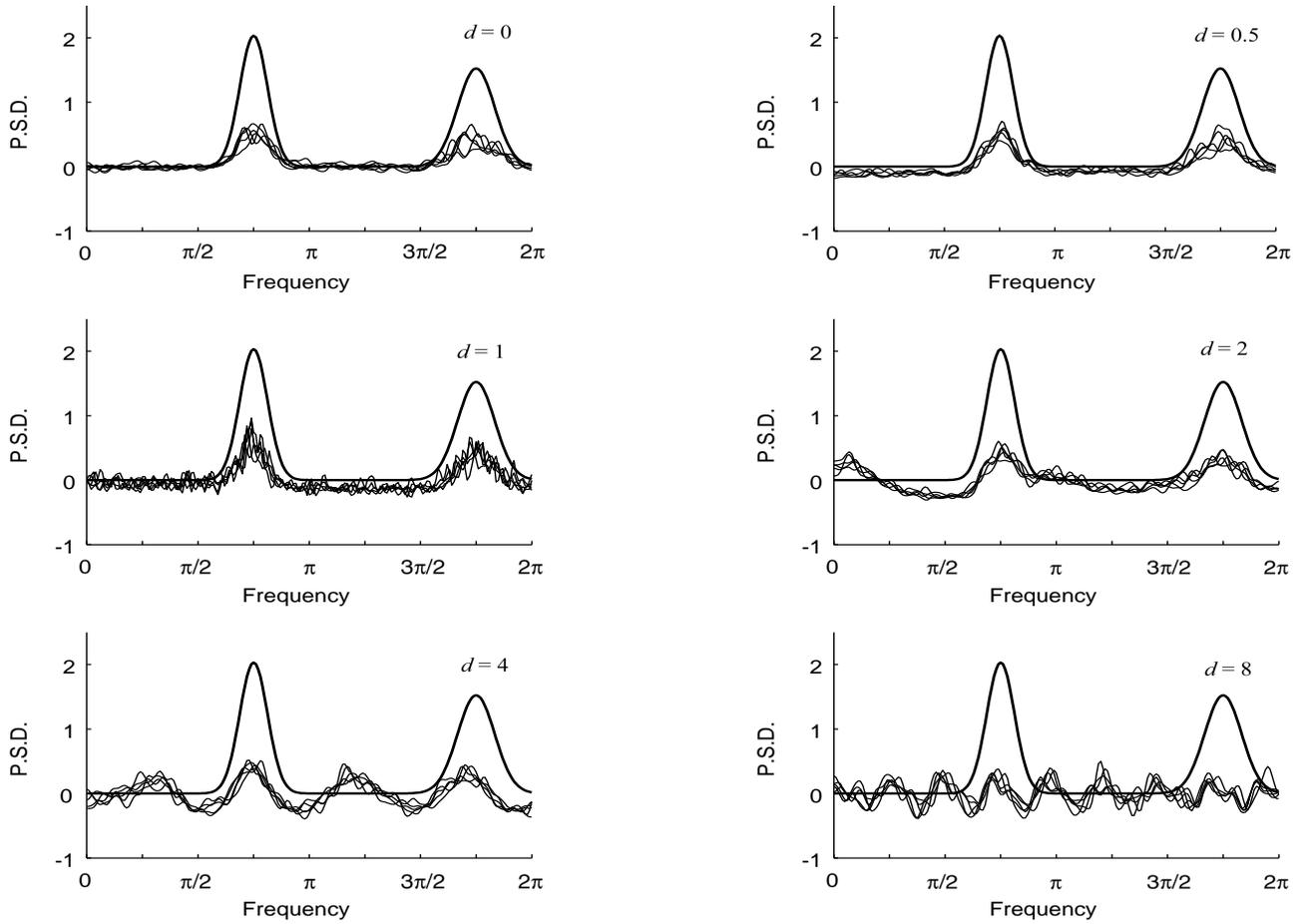}
\caption{Estimates of the power spectral density for $\theta=1$ and
different values of $d$. The bold line represents the true power
spectral density, while the thinner lines represent five typical
estimates.}
\end{figure*}

We consider a continuous time stationary stochastic process $X$ with
mean $0$ and power spectral density $\phi(\cdot)$ given by
$$\phi(\lambda)=\left\{
\begin{array}{l@{\hskip100pt}l}
\multicolumn{2}{l}{\frac{2\sqrt{2}}{\pi}
\left(4e^{-8(4\lambda-3\pi)^{2}/\pi}+3e^{-9(4\lambda-7\pi)^{2}/2\pi}\right)}\\
&\mbox{if $-2\pi\le \lambda\le 2\pi$},\\
0&\mbox{otherwise.}\\\end{array}\right.$$

We consider sampling with a stationary renewal process $\tau$ whose
inter-sample spacing is distributed as $d+T$, where the random
variable $T$ has the exponential distribution with mean~$\theta$.
Note that for $d=0$, this sampling scheme reduces to Poisson
sampling. We assume that $n$ consecutive samples, denoted by
$X(t_1),X(t_2),\ldots,X(t_n)$ are available for estimation.

An estimator of the power spectral density based on the above data,
which is consistent in the special case $d=0$, is given as follows.
\begin{equation}
\begin{split}
\widehat{\phi}_n(\lambda)=&\frac{1}{\pi\beta
n}\sum_{m=1}^{n-1}\sum_{k=1}^{n-m}X(t_k)X(t_{k+m})\\&\hskip20pt
w\left(b_n(t_{k+m}-t_k)\right)\cos\left(\lambda(t_{k+m}-t_k)\right),
\end{split}\label{estim}
\end{equation}
where $w(\cdot)$ is a covariance averaging kernel, and $b_n$ is the
bandwidth of the kernel. Note that this is the estimator proposed by
Masry \cite{Masrynew}, which is an empirical version of the solution
to (\ref{integral}) if one assumes $d=0$ (that is, disregards the
constraint {\bf A4}).

We study the performance of the estimator, $\widehat{\phi}_n(\cdot)$
for the choices
\begin{eqnarray*}
n&=&1000,\\
b_n&=&1/50\\
\mbox{and }w(x)&=&\begin{cases} \frac{1}{2}\left\{1+\cos(\pi
x)\right\}&\mbox{if $-1\le x\le1$,}\\0&\mbox{otherwise,}\end{cases}
\end{eqnarray*}
under the constraint {\bf A4}.

\begin{figure*}[!t]
\noindent\mbox{}\hskip 30pt
\includegraphics[width=6in,height=3in]{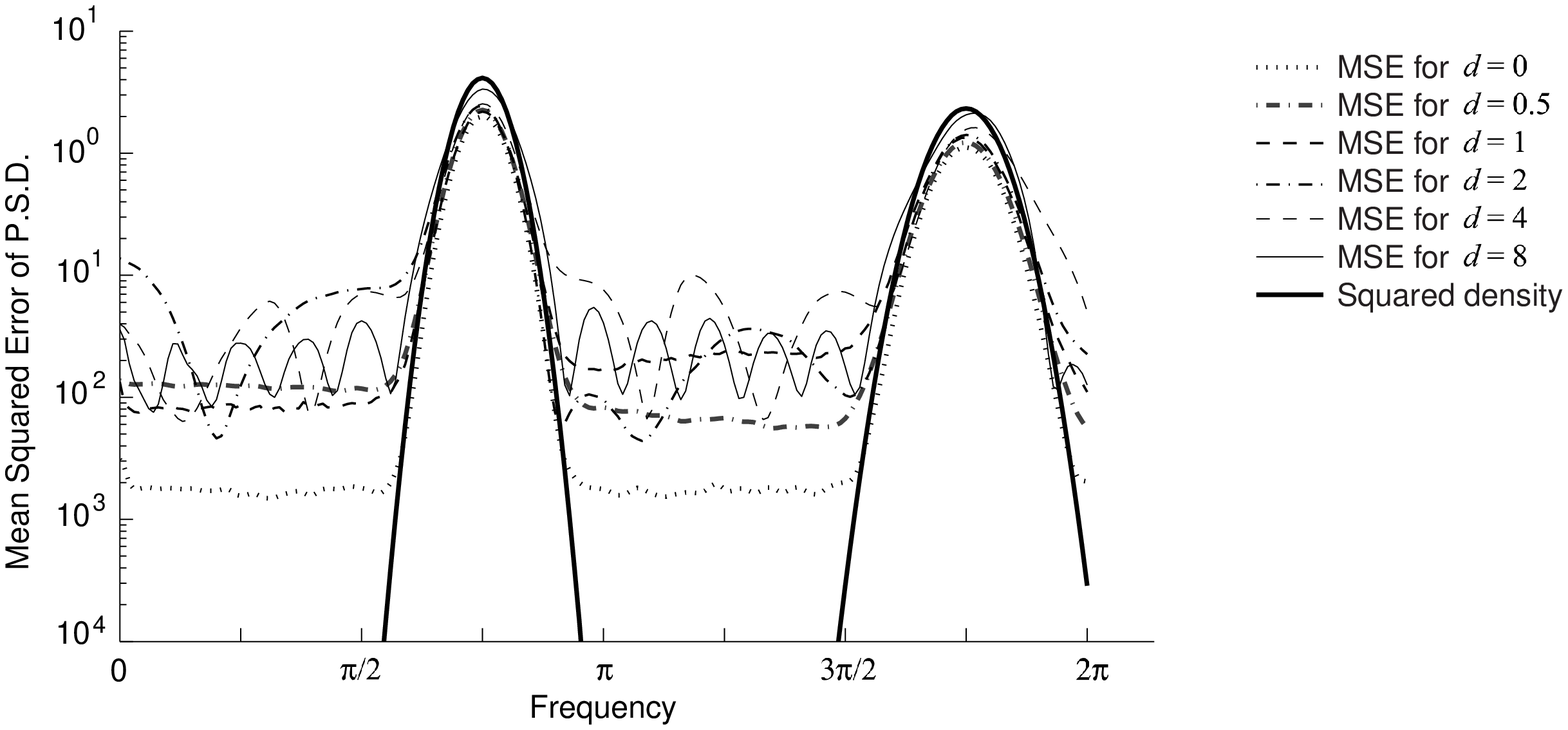}
\caption{Plot of the true squared power spectral density and the
mean squared errors of spectrum estimates (in log scale) based on
500 simulation runs for $\theta=1$ and different values of $d$.}
\end{figure*}

We first investigate how this estimator performs when $d>0$. We
conduct multiple simulation runs for each of the choices $d=0$, 0.5,
1, 2, 4 and 8, together with $\theta=1$. Figure~3 shows spectrum
estimates from five typical simulation runs, along with the true
power spectral density. The plots show how the estimator begins to
perform poorly as one moves away from $d=0$. For larger values of
$d$, the inter-sample spacing is dominated by the constant part.
Therefore, the sampled data resemble that from uniform sampling,
which have the problem of aliasing. As a result, for larger values
of $d$, spurious peaks in greater numbers begin to show in the
estimates. The estimator also assumes larger negative values when
$d$ is larger.

Figure~4 shows the mean squared error (in log-scale) of the estimate
computed in each of the above cases from 500 simulated runs, along
with the squared power spectral density. It is clear that the mean
squared error around the peaks of the power spectral density are of
the same order as the squared power spectral density, and the mean
squared error around the valleys are much larger for $d>0$ than for
$d=0$.

This simulation study indicates that the estimator, which is
consistent in the absence of the constraint on the minimum
inter-sample spacing, can perform poorly in the presence of the
constraint.

The next question we try to answer is: Given the constraint $d=1$
(so that uniform sampling at any sampling rate would necessarily
lead to aliasing), is there an appropriate choice of $\theta$ that
would produce a reasonable estimate of the power spectral density?
In order to answer this question, we again run multiple simulations
for $\theta=0$, 0.05, 0.1, 0.2, 0.5, 1, 2, 5, 10 and 20. In
Figure~5, we present spectrum estimates from five typical simulation
runs in each of these cases, together with the true power spectral
density. For $\theta=0$, i.e., the case of uniform sampling at
sub-Nyquist rate, there is clear evidence of spurious peaks in the
spectrum estimates. A similar occurrence is observed for small
positive values of $\theta$. On the other hand, large values of
$\theta$ give rise to large variability in the estimates.

\begin{figure*}
\noindent
\includegraphics[width=7in,height=9.5in]{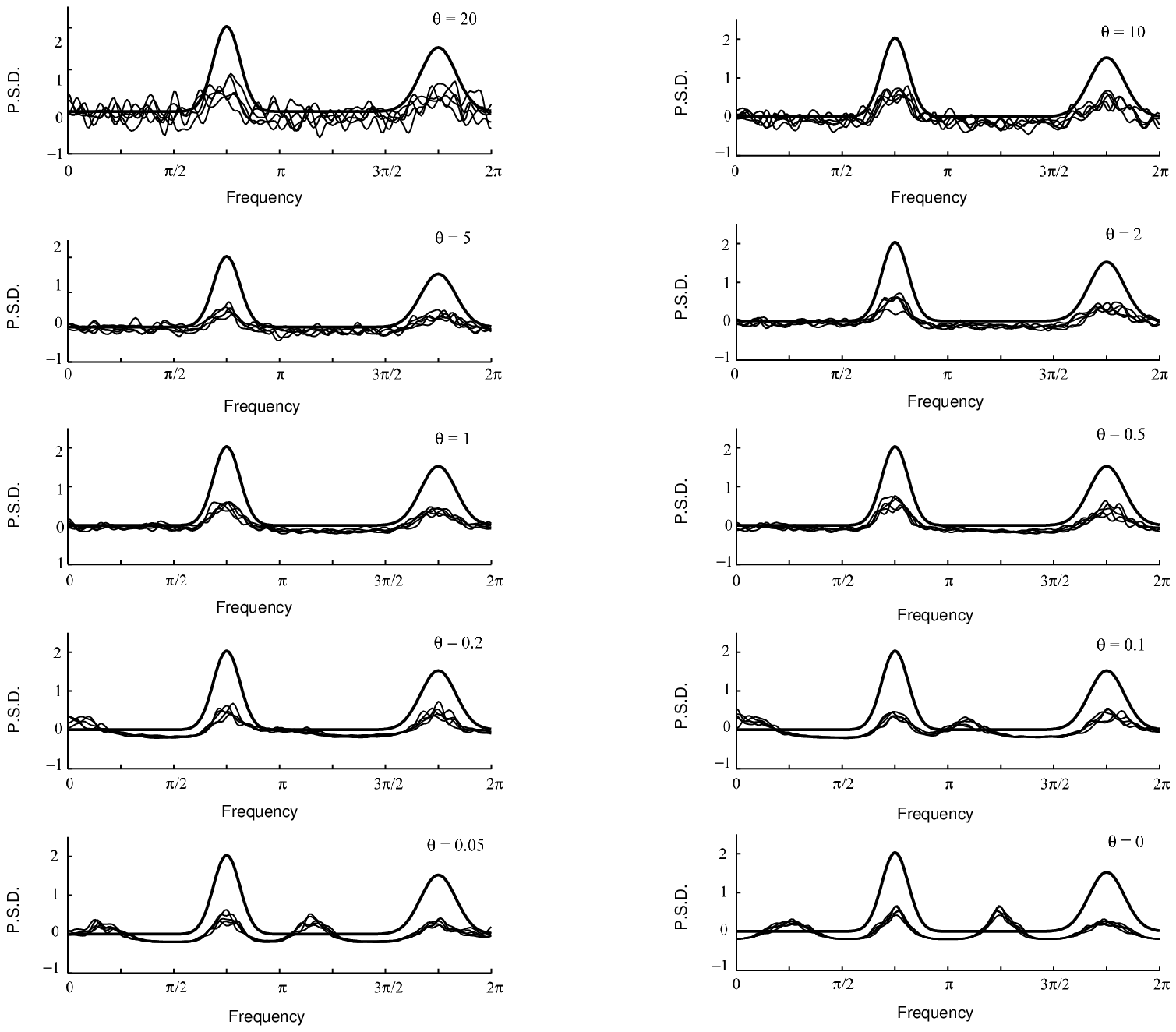}
\caption{Estimates of the power spectral density for $d=1$ and
different values of $\theta$. The bold line represents the true
power spectral density, while the thinner lines represent five
typical estimates.}
\end{figure*}

Figure 6 shows mean squared errors (in log scale) of the estimates
computed in each of the above cases from 500 simulated runs, along
with the squared power spectral density. It transpires that
irrespective of the trade-off between bias and variance observed in
Figure~5, the mean squared errors in all the cases are comparable.
The mean squared error is of the order of the squared value of the
true power spectral density around the peaks, and several orders of
magnitude larger around the valleys.

These findings indicate that the estimator (\ref{estim}) does not
perform well for any choice of $\theta$.

\begin{figure*}[!t]
\noindent\mbox{}\hskip 30pt
\includegraphics[width=6in,height=3in]{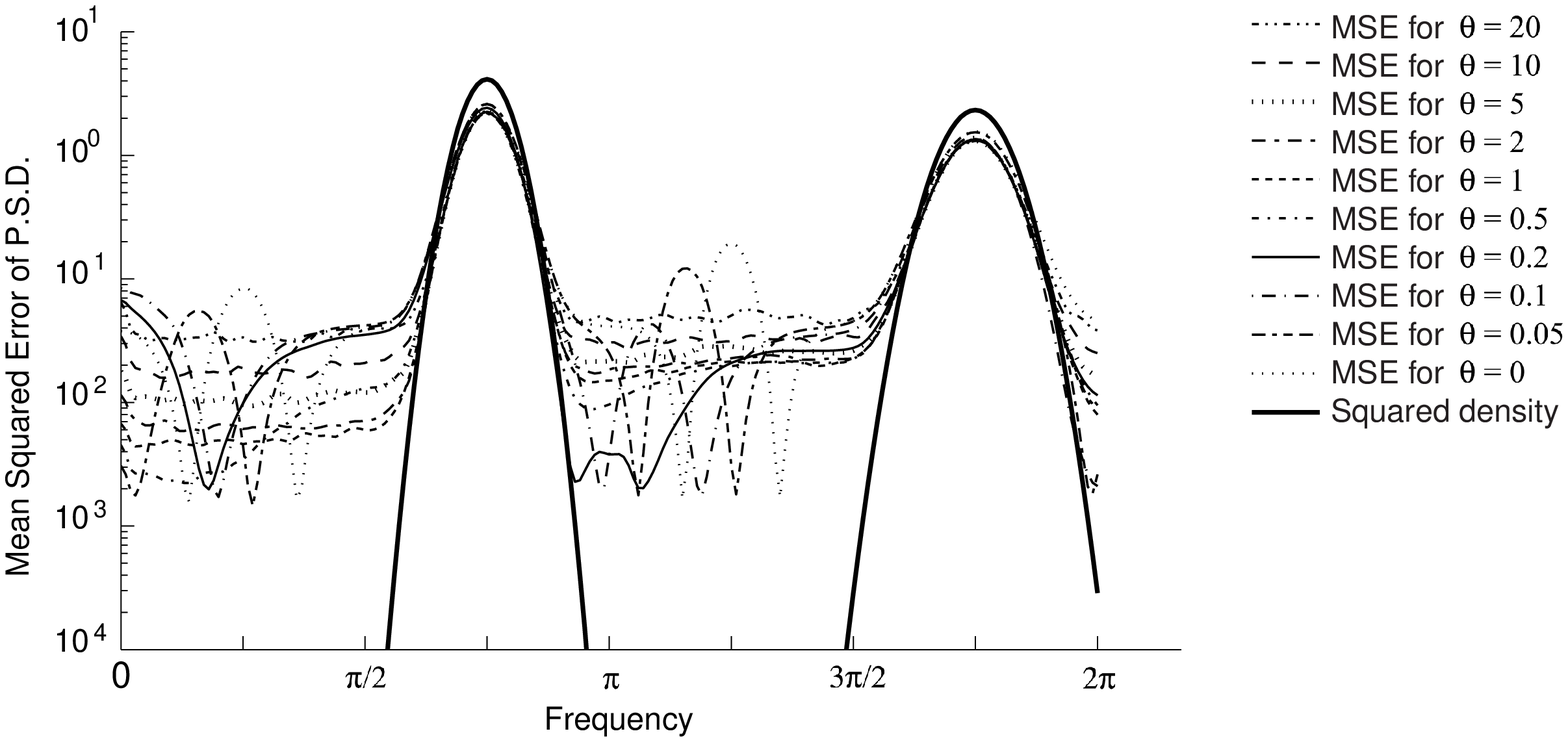}
\caption{Plot of the true squared power spectral density and the
mean squared errors of spectrum estimates (in log scale) based on
500 simulation runs for $d=1$ and different values of $\theta$.}
\end{figure*}

\section{Concluding remarks}\label{S6}

The constraint of a specified minimum inter-sample spacing is a
natural one, in view of technological and economic constraints. We
have come across some interesting findings after formally
incorporating this constraint in the study of aliasing in the
context of spectrum estimation. The most important finding is that
under this constraint, no point process sampling scheme is
alias-free for the class of all spectra -- according to any
definition. It should be noted that the possibility of alias-free
sampling, leading to consistent estimation of the power spectral
density, has traditionally been a major argument in favour of point
process sampling (in contrast with uniform sampling). This argument
does not hold at all in the presence of the above constraint.

We have shown in Section~\ref{S3} that when the inter-sample spacing
is constrained to be larger than a threshold, renewal process
sampling schemes are alias-free for suitably bandlimited spectra
according to Definition 1. The range of bandwidths for alias-free
renewal process sampling for some inter-sample spacing distributions
is smaller than the corresponding range for regular sampling, while
it is larger for some other distributions. On the other hand,
according to Definition 2, all renewal process samplings schemes
satisfying the conditions of Theorem~5 are alias-free for the class
of spectra limited to any finite band.

Masry \cite{Masrynew} pointed out that the plug-in estimator
suggested by Beutler \cite{Beutler} is not provably consistent, and
provided another estimator which is consistent under certain
conditions. The discussion of Section~\ref{S4} shows that these
conditions do not hold under the constraint considered in this paper
-- even when the power spectrum is known to be bandlimited. This
brings us back to square one as far as estimation of power spectral
density is concerned. Since the standard methods would not work well
in the presence of the constraint, as demonstrated through the
simulations of Section~\ref{S5}, there is ample scope for further
research in the area of estimation.

\appendix\setcounter{section}{1}

\noindent {\bf Proof of Theorem 1.} Consider independent, zero mean,
mean square continuous stationary stochastic processes $X_1$, $X_2$
and $X_3$, having covariance functions $C_i(\cdot)$, $i=1,2,3$,
respectively, such that $C_2(0)=C_3(0)$ and $X_2$ and $X_3$ have
different spectral densities belonging to the class $\mathbb{A}$
defined in (\ref{classA}). Consider a sampling point process
$\tau=\{t_n,~n=\ldots,-2,-1,0,1,2,\ldots\}$ satisfying the
Assumptions {\bf A1}, {\bf A2} and {\bf A4}. Let the processes
$X_1+X_2$ and $X_1+X_3$ have spectral distributions
$\Phi_{12}(\cdot)$ and $\Phi_{13}(\cdot)$, respectively, and
covariance sequences of sampled processes
$R_{12}=\{r_{12}(n),~n=\ldots,-2,-1,0,1,2,\ldots\}$ and
$R_{13}=\{r_{13}(n),~n=\ldots,-2,-1,0,1,2,\ldots\}$, respectively.
We have
$$r_{12}(0)=C_1(0)+C_2(0)=C_1(0)+C_3(0)=r_{13}(0).$$

For arbitrary integers $m$ and $n$, let $F_n(x)$ be the distribution
function of $(t_{m+n}-t_m)$. Assumption {\bf A4} implies that
$F_n(x)$ is supported on the interval $\left[|n|d,\infty\right)$. It
follows that, for $n\ne0$,
\begin{equation*}
\begin{split}
r_{12}(n)=&E[X_1(t_{m+n})X_1(t_m)]+E[X_2(t_{m+n})X_2(t_m)]\\
=&\int_{0}^{\infty}C_1(u)dF_n(u)+\int_{0}^{\infty}C_2(u)dF_n(u)\\
=&\int_{|n|d}^{\infty}C_1(u)dF_n(u)+\int_{|n|d}^{\infty}C_2(u)dF_n(u)\\
=&\int_{|n|d}^{\infty}C_1(u)dF_n(u),
\end{split}
\end{equation*}
since $C_2(\cdot)$ is supported on $[-d,d]$. Likewise, $r_{13}(n)$
is also equal to the last expression. This completes the
proof.\mbox{}\hfill$\Box$

\bigskip
\noindent {\bf Proof of Theorem 2.} Consider independent, zero mean,
mean square continuous stationary stochastic processes $X_1$, $X_2$
and $X_3$, having covariance functions $C_i(\cdot)$, $i=1,2,3$,
respectively, such that $C_2(0)=C_3(0)$ and $X_2$ and $X_3$ have
different spectral densities belonging to the class $\mathbb{A}$
defined in (\ref{classA}). Consider a sampling point process
$\tau=\{t_n,~n=\ldots,-2,-1,0,1,2,\ldots\}$ satisfying the
Assumptions {\bf A1}, {\bf A3} and {\bf A4}, and having mean
intensity $\beta$ and reduced covariance measure $\mu_c$. Let the
processes $X_1+X_2$ and $X_1+X_3$ have spectral distributions
$\Phi_{12}(\cdot)$ and $\Phi_{13}(\cdot)$, respectively. As in
Section~\ref{S2}, consider the compound processes
\begin{equation*}
Z_{1j}=\left\{Z_{1j}(B)=\sum_{t_i\in
B}X_1(t_i)+X_j(t_i),~B\in\mathbb{B}\right\},\ j=2,3,
\end{equation*}
which have covariance measures $\mu_{z_{12}}$ and $\mu_{z_{13}}$
given by
\begin{equation*}
\mu_{z_{1j}}(B)=\int_{B}\left\{C_1(u)+C_j(u)\right\}[\beta^2
du+\mu_c(du)],\ j=2,3,
\end{equation*}
respectively.

The reduced covariance measure $\mu_c$ of the point process $\tau$
can be expressed as
$$\mu_c(B)=\beta\delta_0(B)+\beta\int_{B}[dK(u)-\beta
du],~B\in\mathbb{B},$$ %
where
$$K(u)=\sum_{n=1}^{\infty}F_n(|u|),$$
$F_n(u)$ is the conditional probability
$$F_n(u)=\lim_{\epsilon\downarrow 0}P\left[N(t,t+u]\geq n\bigm|N(t-\epsilon,t]\geq 1\right],$$
and $\{N(B),~B\in\mathbb{B}\}$ is the counting process induced by
the process $\tau$ \cite{Beutler1}, \cite{Daley}. Assumption {\bf A4} implies that $K(u)=0$
for $u\in[-d,d]$.

It follows from the above representation of $\mu_c$ that, for each
Borel set $B$, the covariance measures $\mu_{z_{12}}$ is given by
\begin{equation*}
\begin{split}
\mu_{z_{12}}(B)=&\int_{B}C_1(u)[\beta^2 du+\mu_c(du)]\\
&\qquad+\int_{B\cap[-d,d]}C_2(u)[\beta^2 du+\mu_c(du)]
\\=&\int_{B}C_1(u)[\beta^2 du+\mu_c(du)] \\
&\qquad+\beta C_2(0)\delta_0(B\cap[-d,d]).
\end{split}
\end{equation*}
Since $C_2(0)=C_3(0)$, it is clear that the measures $\mu_{z_{12}}$
and $\mu_{z_{13}}$ agree on all Borel sets. This completes the
proof.\mbox{}\hfill$\Box$

\bigskip\noindent
{\bf Proof of Theorem 3.} Here, the sampling process
$\tau=\{t_n,~n=\ldots,-2,-1,0,1,2,\ldots\}$ is such that the
inter-sample spacing $t_{n+1}-t_n$ for different values of $n$ are
independent and identically distributed, say with distribution
function $F(\cdot)$.

Let $\mathbb{S}$ be the class of spectra supported on the closed and finite interval $I$.
Let $X$ be a process as defined in the
theorem, and have the power spectral distribution $\Phi(\cdot)$
belonging to $\mathbb{S}$. The covariance sequence $R$ of the
sampled process is given by
\begin{equation*}
\begin{split}
r(n)=&E[X(t_{m+n})X(t_m)]\\
=&E\left[E\left[X(t_{m+n})X(t_m)\bigm|\tau\right]\right]\\
=&E\left[C(t_{m+n}-t_m)\right]\\
=&E\left[\frac{1}{2\pi}\int_{I}e^{i\lambda(t_{m+n}-t_m)}d\Phi(\lambda)\right]\\
=&\frac{1}{2\pi}\int_{I}E\left(e^{i\lambda(t_{m+n}-t_m)}\right)d\Phi(\lambda).
\end{split}
\end{equation*}
The interchange of the integration is possible by Fubini's theorem,
since the power spectral distribution $\Phi(\cdot)$ and the
probability distribution of $t_{m+n}-t_m$ are both finite. Since the
latter distribution is the $n$-fold convolution of the inter-sample
spacing distribution, we have
\begin{equation}
r(n)=\frac{1}{2\pi}\int_{I}\left[f'(\lambda)\right]^{n}d\Phi(\lambda),
\end{equation}
where $f'(\cdot)$ is the characteristic function of inter-sample
spacing distribution, i.e.,
$$f'(\lambda)=\int_{0}^{\infty}e^{i\lambda y}dF(y),\ \ -\infty<\lambda<\infty.$$

The sampling scheme $\tau$ is alias-free relative to the class of
spectra $\mathbb{S}$ according to Definition 1, if no two different
spectra $\Phi_1$ and $\Phi_2$ belonging to $\mathbb{S}$ produce the
same covariance sequence $R$. Since the sequence $R$ satisfies
$r(-n)=r(n)$, the foregoing condition is equivalent to the
statement:
\begin{equation}\label{criteria}
\begin{split}
&\int_{I}[f'(\lambda)]^{n}(d\Phi_1(\lambda)
-d\Phi_2(\lambda))=0\mbox{ for }n=0,1,2,\ldots\\
& \mbox{ implies that }\Phi_1(\cdot)=\Phi_2(\cdot).
\end{split}
\end{equation}
The above integral with respect to the real variable $\lambda$ can
be written as a complex integral over the contour
\begin{equation}\label{omg}
\Omega=\left\{z:~z=f'(\lambda),~~\lambda\in I\right\}.
\end{equation}
Thus, we can conclude that the sampling scheme $\tau$ is alias free
relative to the class of spectra $\mathbb{S}$ according to
Definition~1 if and only if
\begin{equation}
\begin{split}
&\mbox{For any signed measure $\nu$ defined on the}\\
&\mbox{Borel $\sigma$-field on $\Omega$,}\\
&\int_{\Omega}z^n \nu(dz)=0\mbox{ for }n=0,1,2,\ldots
~\Longrightarrow~ \nu=0.
\end{split}\label{eqv0}
\end{equation}

Note that, since $\Omega$ is the image of the continuous function
$f'(\lambda)$ on the closed and finite interval $I$, the contour
$\Omega$ is compact. Let $C(\Omega)$ be the Banach space of all
complex-valued continuous functions on $\Omega$ equipped with the
supremum norm. Let $\mathbb{M}$ be the set of all signed measures
defined on the Borel $\sigma$-field on $\Omega$. For any
$\nu\in\mathbb{M}$, define the complex valued bounded linear
functional $L_\nu$ defined on $C(\Omega)$ as
\begin{equation}\label{RR}
L_\nu(g)=\int_\Omega g(z)\nu(dz)\quad\mbox{for all } g\in
C(\Omega).
\end{equation}
In terms of these notations, we rewrite (\ref{eqv0}) as
\begin{equation}
\mbox{for any $\nu\in\mathbb{M}$, } ``L_\nu(z^n)=0\mbox{ for
}n=0,1,2,\ldots" ~\Longrightarrow~ \nu=0. \label{eqv1}
\end{equation}
By the Riesz representation theorem, {\it every} bounded linear
functional $L$ on $C(\Omega)$ can be represented as
\begin{equation}\label{RR2}
L(g)=\int_\Omega g(z)\nu_1(dz)+i\int_\Omega g(z)\nu_2(dz)
\quad\mbox{for all } g\in C(\Omega),
\end{equation}
for a unique pair of measures $\nu_1$ and $\nu_2$ in $\mathbb{M}$
\cite{rudin}. It follows that the necessary and sufficient condition
(\ref{eqv1}) is equivalent to the condition:
\begin{equation}
\begin{split}
&\mbox{For any bounded linear functional $L$ on $C(\Omega)$,}\\
&``L(z^n)=0\mbox{ for }n=0,1,2,\ldots"~\Longrightarrow~ L=0.
\end{split}
\label{eqv2}
\end{equation}
The above condition is a statement about the sequence
$\{1,z,z^2,\ldots,\}$ in relation to the Banach space $C(\Omega)$
\{p. 257 of \cite{Davis}\}. By Theorem 11.1.7 of \cite{Davis},
(\ref{eqv2}) is equivalent to the condition:
\begin{equation}
\begin{split}
&\mbox{``The linear span of the sequence $\{1,z,z^2,\ldots,\}$}\\
&\mbox{ is dense in $C(\Omega)$."}
\end{split}\label{eqv3}
\end{equation}
The above condition can be rephrased as: ``Any $g\in C(\Omega)$ can
be expanded in a uniformly convergent sequence of polynomials." By a
result of \cite{mergelyan} (see also \cite{Lavren}), we get the
further equivalent condition:
\begin{equation}
\mbox{``The set $\Omega$ is nowhere dense and does not divide the
plane."}\label{eqv4}
\end{equation}

Since the set $\Omega$ is a curve in the complex plane, it is always
a nowhere dense set. This completes the proof. \mbox{}\hfill$\Box$

\bigskip\noindent
{\bf Proof of Theorem 4.} Let $I=[-1.1\pi/d,1.1\pi/d]$. Consider the
two-point discrete distribution $F$, given by
\begin{equation}
F(t)=\begin{cases} 0&\mbox{for $t<d$,}\\ 0.68&\mbox{for $1\le
t<2.1d$,}\\1&\mbox{for $t\ge 2.1$.}\\
\end{cases}
\label{counter}%
\end{equation}
It follows that the average inter-sample spacing is $1.352d$. Also,
$$f'(\lambda)=0.68e^{i\lambda d}+0.32e^{i2.1\lambda d}.$$
The plot of the imaginary part of $f'(\lambda)$ against the real
part, for $\lambda\in I$, is given in Figure~2. It can be
verified that the graph does not divide the complex plane. The
result follows from Theorem~3.\mbox{}\hfill$\Box$

\bigskip\noindent
{\bf Proof of Theorem 5.} Here, the sampling process
$\tau=\{t_n,~n=\ldots,-2,-1,0,1,2,\ldots\}$ is such that the
inter-sample spacing $t_{n+1}-t_n$ for different values of $n$ are
independent and identically distributed having probability density
function $f(\cdot)$. Let $\beta$ and $\mu_c$ be the mean intensity
and the reduced covariance measure, respectively, of the process
$\tau$. The measure $\mu_c$ can be expressed as
\begin{equation}
\mu_c(B)=\beta\delta_0(B)+\int_{B}\beta[h(u)-\beta]du~\mbox{for
each}~ B\in\mathbb{B},
\label{muc}
\end{equation}
where $h(u)$ is the renewal density function, i.e,
$$h(u)=\sum_{n=1}^{\infty}f^{(n)}(|u|).$$
Note that Assumption {\bf A4} implies that $f(\cdot)$ is supported
on $[d,\infty)$, and so $h(\cdot)$ is supported on
$(-\infty,d]\cup[d,\infty)$. Let us assume, without loss of
generality, that $f(u)>0$ for $u\ge ld$ for some $l\ge1$.

Let $\mathbb{S}$ be the class of bandlimited spectra supported on
$[-\lambda_0,\lambda_0]$. If the sampling scheme $\tau$ is not
alias-free relative to the class of spectra $\mathbb{S}$, then there
exist two zero mean, mean square continuous stationary stochastic
processes $X_1$ and $X_2$ with different power spectral
distributions $\Phi_1(\cdot)$ and $\Phi_2(\cdot)$ such that compound
processes
\begin{equation*}
Z_j=\left\{Z_j(B)=\sum_{t_i\in
B}X_j(t_i),~B\in\mathbb{B}\right\},~j=1,2,
\end{equation*}
have the covariance measures $\mu_{z_1}$ and $\mu_{z_2}$,
respectively, satisfying $\mu_{z_1}=\mu_{z_2}$. Here, for
$B\in\mathbb{B}$, the covariance measures
are given by (see (\ref{muz}) and (\ref{muc}))
\begin{equation*}
\mu_{z_j}(B)=\beta
C_j(0)\delta_0(B)+\beta\int_{B}C_j(u)h(u)du,~j=1,2,
\end{equation*}
where $C_1(\cdot)$ and $C_2(\cdot)$ are the covariance functions of
the processes $X_1$ and $X_2$ respectively. In order that the
covariance measures $\mu_{z_1}$ and $\mu_{z_2}$ are the same, the
point masses at zero, as well as the absolutely continuous parts,
must agree. The equality of the point masses requires
\begin{equation}
C_1(0)=C_2(0).
\label{degpart}
\end{equation}
On the other hand, equality of the absolutely continuous parts means
$$C_1(u)h(u)=C_2(u)h(u)\mbox{ for }-\infty<u<\infty.$$
Since $h(u)>0$ for the $|u|\ge ld$, we have
\begin{equation}
C_1(u)=C_2(u)~\mbox{for}~|u|\ge ld. \label{2cov}
\end{equation}

If the processes $X_1$ and $X_2$ have spectra limited to the band
$[-\lambda_0,\lambda_0]$, then the covariance function $C_j(\cdot)$
for $j=1,2$ can be expressed as \cite{OpSch}
\begin{equation}
C_j(u)=\frac{1}{T}\sum_{n=-\infty}^{\infty}C_j(nT)\mbox{sinc}\left(\frac\pi
T (u-nT)\right), \label{repr}
\end{equation}
where $T=\frac{2\pi}{2\lambda_0}$ and
$$\mbox{sinc}(x)=\begin{cases}
\frac{\sin x} x&\mbox{if $x\ne0$,}\\ 1&\mbox{if
$x=0$.}\\\end{cases}$$

Let $k=[ld/T]$, where $[u]$ represents the integer part of the real
number $u$. It follows from (\ref{degpart})--(\ref{repr}) that
\begin{equation*}
\begin{split}
&\hskip-30pt
C_1(u)-C_2(u)\\=&\sum_{n=-k}^{k}\left\{C_1(nT)-C_2(nT)\right\}
\mbox{sinc}\left(\frac\pi T (u-nT)\right)\\
&+\sum_{|n|>k}\left\{C_1(nT)-C_2(nT)\right\}\mbox{sinc}\left(\frac\pi T (u-nT)\right)\\
=&\sum_{n=1}^{k}\{C_1(nT)-C_2(nT)\}\\
&\times\left(\mbox{sinc}\left(\frac\pi T (u-nT)\right) %
+\mbox{sinc}\left(\frac\pi T (u+nT)\right)\right).
\end{split}
\end{equation*}
By using the fact that $\sin(k\pi+\theta)=(-1)^{k}\sin\theta$ for
all integer $k$, we have for $\alpha=u/T-[u/T]>0$,
\begin{equation*}
\begin{split}
&\hskip-7pt C_1(u)-C_2(u)\\
=&\sum_{n=1}^{k}\{C_1(nT)-C_2(nT)\}\\
&\times\!\left(\frac{\sin\left\{\left(-n+\!\left[\frac u
T\right]\right)\pi+\alpha\pi\right\}}{\frac{\pi}{T}(u-nT)}
      +\frac{\sin\left\{\left(n+\!\left[\frac u
T\right]\right)\pi+\alpha\pi\right\}}{\frac{\pi}{T}(u+nT)}\!\right)\\
=&(-1)^{[u/T]}\sin(\alpha\pi)\\
&\times\sum_{n=1}^{k}\{C_1(nT)-C_2(nT)\}
\left(\frac{(-1)^{-n}}{\frac{\pi}{T}(u-nT)}+\frac{(-1)^{n}}{\frac{\pi}{T}(u+nT)}\right).
\end{split}
\end{equation*}
Since $(-1)^{n}=(-1)^{-n}$ for each integer $n$, we have
\begin{equation*}
\begin{split}
&\hskip-30pt C_1(u)-C_2(u)\\
=&(-1)^{[u/T]}\frac{2uT}{\pi}\sin(\alpha\pi)\\
&\times\sum_{n=1}^{k}[(-1)^{n}\{C_1(nT)-C_2(nT)\}]\frac{1}{u^2-n^{2}T^{2}}.
\end{split}
\end{equation*}
Let $v_n=[(-1)^{n}\{C_1(nT)-C_2(nT)\}]$. In view of (\ref{2cov}),
the above equation implies that
\begin{equation}
\sum_{n=1}^{k}\frac{v_n}{u^2-n^{2}T^{2}}=0,\label{ratiopol}
\end{equation}
for
$u\in\{(ld,(k+1)T)\}\cup\left\{\cup_{m=k+1}^\infty(mT,(m+1)T)\right\}$.

Note that the function on the left hand side of (\ref{ratiopol}) is
a ratio of polynomials. The polynomial in the numerator has degree
$2k-2$, while the denominator is bounded over the domain of the
function. Thus, the ratio of the polynomials can be zero at most at
$2k-2$ points. Therefore, the fact that this function assumes the
value 0 everywhere on the interval\linebreak $((k+1)T,(k+2)T)$
implies that the polynomial in the numerator is identically equal to
zero. Thus, the ratio of the polynomials is identically zero. Hence,
$$\sum_{n=1}^{k}\frac{v_n}{u^2-n^{2}T^{2}}=0,\mbox{ for }
u\in\bigcup_{m=0}^\infty(mT,(m+1)T).$$ %
By considering the limit of the left hand side as $u\downarrow nT$,
it is found that $v_n=0$ for $n=1,\ldots,k$, that is,
$$C_1(nT)=C_2(nT)\mbox{ for }|n|=1,\ldots, k.$$
According to (\ref{degpart}), the above equality holds for $n=0$,
while (\ref{repr}) and (\ref{2cov}) imply that it holds for
$|n|=k+1,k+2,\ldots$. Thus, $C_1(nT)=C_2(nT)$ for all $n$. It
follows from (\ref{repr}) that $C_1(u)=C_2(u)$ for each $u$, which
contradicts the assumption that $C_1$ and $C_2$ are different. So
the sampling scheme $\tau$ is alias-free for the class of the
spectra $\mathbb{S}$.\hfill$\Box$

\section*{Acknowledgements}
The authors gratefully acknowledge suggestions and technical help
from Professor B.V. Rao of the Chennai Mathematical Institute.

\end{document}